\documentclass[11pt]{amsart}
\usepackage{amsmath,amsfonts,amssymb,amsxtra,latexsym,amscd,enumerate,amsthm}

\usepackage{latexsym}
\usepackage{epsfig}
\usepackage{verbatim}

\def\a{\alpha}

\def\p{\Phi}
\def\v{\varphi}

\def\l{\lambda}
\def\o{\omega}

\def\R{\mathbb{R}}

\def\P{\mathbb{P}}
\def\p{\mathcal{P}}

\def\I{\mathcal{I}}

\def\A{\mathcal{A}}
\def\B{\mathcal{B}}
\def\TT{\mathbb{T}}
\def\N{\mathbb{N}}

\def\f{\mathcal{F}}
\def\W{\mathcal{W}}

\def\beq{\begin{equation}}
\def\eeq{\end{equation}}

\def\beq{\begin{equation}}
\def\eeq{\end{equation}}

\newtheorem{t1}{Theorem}
\newtheorem{l1}{Lemma}
\newtheorem{p1}{Proposition}

\begin{document}
\title[]{On the pointwise convergence of the sequence of partial Fourier Sums
along lacunary subsequences}

\author{Victor Lie}

\date{\today}
\address{Department of Mathematics, Princeton, NJ 08544-1000 USA}

\email{vlie@math.princeton.edu}

\address{Institute of Mathematics of the Romanian Academy, Bucharest, RO
70700 \newline \indent  P.O. Box 1-764}


\maketitle

\begin{abstract}
In his 2006 ICM invited address, Konyagin mentioned the following conjecture: if $S_n f$ stands for the
$n$-th partial Fourier sum of $f$ and $\{n_j\}_j\subset \N$ is a lacunary sequence, then $S_{n_j} f$ is a.e. pointwise convergent for any
$f\in L\log\log L$. In this paper we will show that
$\|\,\sup_{j} |S_{n_j}(f)|\,\|_{1,\infty}\leq C\, \|f\|_{1}\,\log\log (10+\frac{\|f\|_{\infty}}{\|f\|_1})\:.$ As
a direct consequence we obtain that $S_{n_j}f\,\rightarrow\,f $ a.e. for $f\in L\log\log L\log\log\log L$.
The (discrete) Walsh model version of this last fact was proved by Do and Lacey but their methods
do not (re)cover the (continuous) Fourier setting. The key ingredient
for our proof is a tile decomposition of the operator $\sup_{j} |S_{n_j}(f)|$ which depends on both the function $f$ and on the lacunary
structure of the frequencies. This tile decomposition, called $(f,\l)-$lacunary, is directly adapted to the context of our problem, and,
combined with a canonical mass decomposition of the tiles, provides the natural environment to which the methods developed by the author
in ``On the Boundedness of the Carleson Operator near $L^1$" apply.

\end{abstract}
$\newline$

\section{\bf Introduction}

This paper extends the line of research addressed in \cite{lv7} to the problem regarding the pointwise convergence of the \textit{lacunary}\footnote{From now on, we will ``slightly abuse" the terminology and refer to any lacunary subsequence of the sequence of partial Fourier sums simply as a lacunary Fourier series.} Fourier Series near $L^1$. This problem was extensively studied by V. Konyagin (\cite{ko1}, \cite{ko2}) and it came to our attention when reading \cite{LaDo}.

In what follows, we will mention several historical facts about the evolution of the subject treated here. In the extended framework, the problem of the pointwise convergence of the \textit{full}\footnote{Here we refer to the entire sequence of partial Fourier sums.} sequence of the partial Fourier sums evolved as described briefly below: in
\cite{Kol1}, Kolmogorov showed that there exists $f\in L^1(\TT)$ such that $\{S_n(f)\}_n$ is almost everywhere divergent. Next, more than forty years later, L. Carleson (\cite{c1}) gave a positive answer to Lusin's conjecture - \textit{i.e.} - the Fourier series of a function $f\in L^2(\TT)$ is almost everywhere convergent. Then, Hunt (\cite{hu}) extended his result to the setting of the $L^p(\TT)$ spaces for $1<p<\infty$. The Carleson-Hunt theorem was later reproved by C. Fefferman (\cite{f}) and M. Lacey and C. Thiele (\cite{lt3}). The remaining fundamental question
is:`` What can one say about the behavior of the Fourier Series near $L^1$? " or, generally speaking, how should the pointwise convergence story reveal between the negative result of Kolmogorov ($p=1$) and the positive result of Carleson-Hunt ($p>1$)? Several steps were made in the direction of clarifying this story (\cite{sj3},\cite{So1},\cite{So2},\cite{An},\cite{Ar}) with a unifying perspective on these previous results offered by the author in \cite{lv7}. Essentially, at this time, the best result is due to Antonov (\cite{An}), and asserts the pointwise convergence of the Fourier series for functions belonging to the Orlicz space
$L\log L\log\log\log L(\TT)$.

The challenging ``mystery of the story" is represented by
$\newline$

\noindent\textbf{Conjecture 1.} \textit{The sequence of the partial Fourier sums $\{S_n(f)(x)\}_n$ is a.e. $x\in\TT$ convergent for any $f\in L\log L(\TT)$.}
$\newline$

In a different but close in spirit direction one may ask - ``what is the largest Banach function space for which one has
(pointwise a.e.) convergence of the \textit{lacunary} Fourier Series ?" In a symmetric treatment with that of the previous topic, we start by mentioning that Gosselin  proved in \cite{Gos} that for any increasing sequence $\{n_j\}_j\subseteq\N$ there is a function $f\in L^1(\TT)$ such that
$\sup_{j}|S_{n_j}(f,x)|=\infty\:\:\textrm{a.e.}\:x\in\TT\,.$ Surprisingly enough though, if $f\in H^1(\TT)$ and $\{n_j\}_j\subseteq\N$ is lacunary, then we have that $\{S_{n_j}(f)\}_j$ is a.e. pointwise convergent.\footnote{Notice that $H^1(\TT)$ is not a Banach function space.} Moreover, according to a result of Zygmund (\cite{Zyg}), the above conjecture is true if one merely restricts to the convergence of lacunary subsequences of the partial Fourier sums!

Hence, we do expect a significantly better behavior for the pointwise convergence of the lacunary Fourier series. Indeed, one hopes for the following
to be true
$\newline$

\noindent$\textbf{Conjecture 2.}$ (Konyagin,\cite{ko1}) \textit{Let $\{n_j\}_j\subset\N$ be a lacunary sequence. Then there exists
$C>0$ absolute constant such that if $f\in L\log\log L(\TT)$ the following holds}
\beq\label{conjKonyag}
\|\,\sup_{j} |S_{n_j}(f)|\,\|_{1,\infty}\leq C\, \|f\|_{L\,\log\log L}\:.
\eeq
\textit{As a consequence one also has}
$$S_{n_j}(f,x)\:\stackrel{j\rightarrow\infty}{\longrightarrow}\:f(x)\:\:\:\:\:\textrm{a.e.}\:\:x\in\TT\:.$$

It is worth mentioning that, if true, Conjecture 2 is sharp since in \cite{ko2} Konyagin proves that for any increasing sequence $\{n_j\}_j\subset\N$ and any increasing function $\phi:\R_{+}\rightarrow\R_{+}$ such that $\phi=o(u\,\log\log u)$ when $u\rightarrow\infty$ there
is a function $f\in\phi(L)$ such that $\sup_{j}|S_{n_j} (f,x)|=\infty$ for all $x\in\TT$.

The main result of this paper is
$\newline$

\noindent\textbf{Main Theorem.} \textit{Let  $f\in L^{\infty} (\TT)$  and $\{n_j\}_j\subset\N$ a lacunary sequence. Then we have
\beq\label{Konyagin}
\|\,\sup_{j} |S_{n_j}(f)|\,\|_{1,\infty}\leq C\, \|f\|_{1}\,\log\log (10+\frac{\|f\|_{\infty}}{\|f\|_1})\:,
\eeq
where here $C>0$ does not depend on $f$.
As a consequence:}

\textit{i) Conjecture 2 is true for any $f=\chi_F$ where $F\subseteq\TT$ measurable.}

\textit{ii) In general we have that
\beq\label{K1}
\|\,\sup_{j} |S_{n_j}(f)|\,\|_{1,\infty}\lesssim \|f\|_{L\log\log L\log\log\log L}\:,
\eeq
and hence
$$S_{n_j}(f,x)\:\stackrel{n\rightarrow\infty}{\longrightarrow}\:f(x)\:\:\:\:\:\textrm{a.e.}\:\:x\in\TT$$
for any $f\in L\log\log L\log\log\log L$.}
$\newline$

Our result shows that Konyagin's conjecture is true up to a $\log\log\log$ factor.
If one modifies Conjecture 2 addressing it in the Walsh-Fourier setting, then \eqref{K1} was shown to hold in \cite{LaDo}.
However, the methods used by Do and Lacey are reducing the problem to a projection argument that does not seem to extend to the continuous case treated by us. That is why, we will embrace a different path developing a tile discretization adapted to the nature of our problem. For more details on the antithesis between \cite{LaDo} and the present paper see Remarks section.
$\newline$

\section{Discretization of the operator}

Since the maximal operator under discussion is nothing else than a \textit{lacunary} version of the Carleson operator,
as usual in such context, we will use time-frequency methods to analyze it.

Now, the study of our operator
\beq\label{carlac}
S_{lac}f(x):=\sup_{j}|S_{n_j}f(x)|\:\:\:\:\:\:\:f\in C^1(\TT)
\eeq
may be canonically reduced to the analysis of
\beq\label{carlac1}
Tf(x):=\sup_{j\in\N}\,|\int_{\TT} \frac{1}{x-y}\,e^{i\,n_j\,(x-y)}\,f(y)\,dy| \:,
\eeq
where here $\{n_j\}_j$ is a prescribed lacunary sequence of positive integers.

Applying Fefferman's approach (\cite{f}) we perform the following steps:
\begin{itemize}
\item linearize our operator and write\footnote{For technical reasons we will erase the term $N(x)\,x$ in the phase of the exponential, as later
in the proof this will simplify the structure of the adjoint operators $T_P^{*}$.}
$$Tf(x):=\int_{\TT} \frac{1}{x-y}\,e^{-i\,N(x)\,y}\,f(y)\,dy\,,$$
where here $N:\:\TT\rightarrow\:\{n_j\}_j\;$ measurable function.

\item use the dilation symmetry of the kernel and express
$$\frac{1}{y}=\sum_{k\geq 0} \psi_k(y)\:\:\:\:\:\:\:\:\:\forall\:\:0<|y|<1\:,$$
where  $\psi_k(y):=2^{k}\psi(2^{k}y)$ (with $k\in \N$) and $\psi$ an odd
$C^{\infty}$ function such that
$\operatorname{supp}\:\psi\subseteq\left\{y\in \R\:|\:2<|y|<8\right\}$.

\item write $$Tf(x)=\sum_{k\geq 0}\int_{\TT}e^{-i\,N(x)\,y}\,\psi_{k}(x-y)\,f(y)\,dy\:.$$

\item partition the time-frequency plane in tiles (rectangles of area one) of the form $P=[\o,I]$ with $\o,\,I$
dyadic intervals\footnote{With respect to the canonical dyadic grids on $\R$ and respectively $\TT$.} such that $|\o|=|I|^{-1}$.
Set the collection of all such tiles as $\P$.

\item to each  $P=[\o,I]\in\P$ we assign the set $E(P):=\left\{x\in I\:|\:N(x)\in
P\right\}$ that is responsible for the ``weight" of the tile - $\frac{|E(P)|}{|I|}$
depending on which we will later realize a first partition of the set $\P$.

\item for $P=[\o,I]\in\P$ with $|I|=2^{-k}$ ($k\geq0$) we define the operators
$$T_{P}f(x)=\left\{\int_{\TT}e^{-i\,N(x)\,y}\,\psi_{k}(x-y)\,f(y)\,dy\right\}\chi_{E(P)}(x)\:,$$
and conclude that
\beq\label{discret}
Tf(x)=\sum_{P\in\P} T_P f(x)\,.
\eeq
\end{itemize}

Notice that if we think to $N:\:\TT\rightarrow\:\{n_j\}_j\;$ as a predefined measurable function then the above decomposition is
\textit{independent} on the function $f$. Using this perspective will be enough to show that the bounds on $T$ do not depend on $N$.

\section{Discretization of the families of tiles $\P$}

In this section we will decompose the family of tiles $\{P=[\o_P, I_P]\}_{P\in\P}$ according to two different concepts:
\begin{itemize}
\item the \textbf{``weight" (or mass)} of a tile - $A(P)=\frac{|E(P)|}{|I_P|}$; this decomposition is thus \textit{independent} on the function $f$.

\item the \textbf{$\o_P$ (frequency) localization  versus $I_P$ (spacial) concentration of $f$}; this decomposition (called \textbf{$(f,\l)-$lacunary}\footnote{The parameter $\l$ quantifies the size of the Hardy-Littlewood maximal function of $f$.}) depends on both the relative
position of each tile $P$ with respect to the real axis\footnote{Here is precisely where the structure of our operator intervenes \textit{i.e.} the supremum is taken only over \textit{lacunary} frequencies.} and on the information carried by the spacial support of $P$ with respect to $f$.
\end{itemize}

\subsection{The mass decomposition}

This decomposition was first developed in the seminal approach of Fefferman (\cite{f}). Here though we will make use of a refinement of it as appearing in \cite{lv7} when treating the exceptional sets in the Polynomial Carleson operator discretization. This refined decomposition is also carefully described in \cite{lv3}. For this reason we will skip the details of this decomposition (the interested reader should consult Section 5 in \cite{lv3}) and only mention that - heuristically - as an output of this procedure we will be able to write\footnote{At a more precise level, each family $\P_n$ can be reduced to a $BMO-$forest of $n^{th}$ generation - again, for the definition see Section 4 in \cite{lv3}.}
$$\P=\bigcup_{n\in\N}\P_n\:,$$
with each $\P_n\cong\{P\in\P\,|\,A(P)\approx 2^{-n}\}\:.$

\begin{comment}
Te iubesc, Tapi-Tapi! Cu Drag, TTH (Los Angeles, 16 martie 2012)(forever)
\end{comment}

\subsection{The $(f,\l)-$lacunary decomposition}

As mentioned at the beginning, in this section we will perform a second decomposition of our tiles depending on the size/localization
of the function\footnote{Thus, notice that this second decomposition is dependent on $f$ hence each function will involve different partitions. However for notational simplicity we will not write explicitly the $f$ dependence in our decomposition.} $f$ and on the geometric location of the tiles with respect to the origin that accounts for the lacunary structure of our maximal operator.

  For expository reasons in what follows we will only refer to the case when $f=\chi_F$ (here $F\subseteq\TT$ measurable). This will give us a ``simplified" picture of our decomposition which still encapsulates the \textit{essence} of the matter. For a general $f$ the required modifications will be discussed in Section 4.2..

\begin{comment}
Let $\l\in (0,1)$ be a fixed parameter and $\I$ be the collection of maximal dyadic intervals $I$ satisfying\footnote{Here we assume wlog that $\frac{|F|}{|G|}<1$ since the other
case reduces trivially to the $L^2$ boundedness of the Carleson operator.}
\beq\label{max}
\frac{|I\cap F|}{|I|}\geq\lambda=\frac{|F|}{|G|}\;.
\eeq
\end{comment}

Let $\l\in (0,1)$ be a fixed parameter.
For each $k\in\N$ set $\I_k$ the collection of maximal dyadic intervals $I$ such that
$$\frac{|F\cap I|}{|I|}> \l\,2^{-k}\:,$$
and set $\bar{\I}_k=\bigcup_{I\in\I_k}\,I$.

Observe that for any $I\in\I_k$ with $10|I|<|\TT|$ there exists $J\in\I_{k+1}$ such that $I\subsetneq J$ and hence $\bar{\I}_k\subsetneq \bar{\I}_{k+1}$.

Passing to the tile discretization algorithm, we first want to isolate (remove) the family of tiles that are not well separated. For this, assuming without loss of generality\footnote{The sequence $\{n_k\}_k$ lacunary implies $\lim \textrm{inf}_{k\rightarrow\infty}\frac{n_{k+1}}{n_k}=\a>1$.} that $n_k= \a^k$ with $\a\in\N,\:\a\geq 2$ we define\footnote{Throughout this paper we will use the following standard notation: if $I$ is an (open) interval having the center $c$, then for any $b>0$ we set $b\,I:=(c-\frac{b\,|I|}{2},\,c+\frac{b\,|I|}{2})$.}
\beq\label{clustertile}
\P_{cluster}:=\{P=[\o_P, I_P]\in\P\,|\,0\in 10\,\a\,\o_{P}\}\:.
\eeq
Next, we set
$$\p_{k,O}:=\{P_O=[\o_{P_O}, I_{P_O}]\in\P\,|\ \exists\,I\in\I_k\:\:s.t.\:\:I_{P_O}=I\:\textrm{and}\:0\in2\o_{P_O}\}\:.$$

In what follows we will split the entire family of tiles $\P_{sep}:=\P\setminus\P_{cluster}$ relative to the structure offered by the sets $\{\p_{k,O}\}_{k\in\N}$.

Since our procedure will involve the support of the adjoint operators $\{T_P\}_{P\in\P_{sep}}$ we first isolate an elementary piece $T_P$ (and the corresponding $T_P^*$) and briefly introduce several notations that we will use in our construction:

For $P=[\o_P,I_P]\in\P$ we set $c(I_P)$ the center of the interval $I_P$ and define $I_{P^*}=[c(I_P)-\frac{17}{2}|I_P|,\,c(I_P)-\frac{3}{2}|I_P|]\cup[c(I_P)+\frac{3}{2}|I_P|,\,c(I_P)+\frac{17}{2}|I_P|]$; we then have the following properties:
\beq\label{support}
\textrm{supp}\,T_{P}\subseteq I_P\:\:\:\textrm{and}\:\:\:\textrm{supp}\,T_{P}^{*}\subseteq I_{P^*}\:.
\eeq
Notice that we can express the set containing the support of $T_{P}^{*}$ as $$I_{P^*}=\bigcup_{r=1}^{14}I_{P*}^r$$ with each $I_{P^*}^r$ a dyadic
interval of length $|I_P|$.

This being said let us resume our tile decomposition; for $P_O\in\p_{k,O}$ we set\footnote{We perform this decomposition for those $k\in\N$ for
which $\bar{\I}_k\subsetneq\TT$ with the obvious modifications required for the case $\bar{\I}_k=\TT$; for the remaining $k'$s, we will set
$\p_k^{2}(P_O)=\p_k^{1}(P_O)=\emptyset$.}
\beq\label{l2fam}
 \p_k^{2}(P_O):=\left\{\begin{array}{cl}P=[\o_{P}, I_{P}]\\P\in\P_{sep}\end{array}\:\big|\:\begin{array}{cl}\exists\:\:I_{P^*}^r\supseteq I_{P_O}\:\textrm{and}\:2\o_{P}\cap2\o_{P_O}=\emptyset\:\\
  \textrm{if}\:I\in\I_{k+1}\:s.t.\:I\cap I_{P^*}^r\not=\emptyset\:\textrm{then}\:I\supseteq I_{P^*}^r\end{array}\right\},
\eeq
 and
\beq\label{l1fam}
\p_k^{1}(P_O):=\left\{\begin{array}{cl}P=[\o_{P}, I_{P}]\\P\in\P_{sep}\end{array}\:\big|\:\begin{array}{cl}\exists\:\:I_{P^*}^{r'}\supseteq I_{P_O}\:\textrm{and}\:2\o_{P}\cap2\o_{P_O}\not=\emptyset\:\\
  \textrm{if}\:I\in\I_{k+1}\:s.t.\:I\cap I_{P^*}^{r'}\not=\emptyset\:\textrm{then}\:I\supseteq I_{P^*}^{r'}\end{array}\right\},
\eeq
where here $r,\,r'\in\{1,\ldots\,14\}$.

Notice that we have now that
 \beq\label{tot}
 \P=\P_0\cup\P_{cluster}\cup\bigcup_{k\in\N}\bigcup_{P_O\in\p_{k,O}}\p_k^{2}(P_O)\cup\p_k^{1}(P_O)\:,
 \eeq
where $\P_0$ is a collection of tiles such that $\forall\:P=[\o_P,\:I_P]\in\P_0$ one must have
\beq\label{rest}
|I_{P^*}\cap F|=0\:\:\:\:\textrm{or}\:\:\:\:I_{P^*}\subset 1000\, F_{bad}\:\: \textrm{with}\:\: F_{bad}:=\{x\,|\,M(\chi_F)(x)> \frac{1}{2}\l\}\:.
\eeq

Also it is worth mentioning that \eqref{tot} does not express $\P$ as a disjoint union (partition) of sets; more precisely it is possible that
 $\p_k^{l}(P_O)\cap\p_{k'}^{l'}(P'_O)\not=\emptyset$ for some $k,\,k'\in\N$, $l,l'\in\{1,2\}$, $P_O\in\p_{k,O}$ and $P'_O\in\p_{k,O}$. On the other hand one should notice that for any $P\in\P$ we have
 $$\#\{k\,|\,P\in \p_k^{l}(P_O)\,\textrm{for some}\:P_O\in \p_{k,O}\,,\:l\in\{1,2\}\}\leq 14\:.$$

However, a key observation is that when transferred in the setting of the initial problem our tile decomposition behaves as good as a partition\footnote{Here is the key point where we are taking advantage that we only need $L^{1,\infty}$ estimates for our operator.}. To see this,
given any measurable set $\bar{G}\in\TT$ choose\footnote{Here we assume wlog that $\frac{|F|}{|\bar{G}|}<1$ since the other
case reduces trivially to the $L^2$ boundedness of the Carleson operator.} $\l\approx\frac{|F|}{|\bar{G}|}$ and define $$G:=\bar{G}\setminus1000\,F_{bad}\:.$$
From the above construction we notice $|\bar{G}|\approx|G|\:.$

With these facts, making use of one more observation:
$$\forall\:I\:\textrm{interval}\:I\subset \bar{\I}_{k+1}\:\textrm{and}\:I \cap\bar {\I}_{k}=\emptyset\:\:\Rightarrow\:\:|I\cap F|=0\:,$$
we conclude
\beq\label{partition}
\begin{array}{rl}
\|\sum_{P\in\P} T_P^{*}(\chi_G)\|_{L^1(F)}&\leq \|\sum_{P\in\P_{cluster}} T_P^{*}(\chi_G)\|_{L^1(F)}\\
&+\,\|\sum_{k}\sum_{I_{P_O}\in\I_k}\chi_{I_{P_O}}{T^{\p_k^{2}(P_O)}}^{*}(\chi_G)\|_{L^1(F)}\\
&+\,\|\sum_{k}\sum_{I_{P_O}\in\I_k}\chi_{I_{P_O}}{T^{\p_k^{1}(P_O)}}^{*}(\chi_G)\|_{L^1(F)}\;.
\end{array}
\eeq

\section{The proof of the Main Theorem}

In this section we prove our main result. We will do this in two steps, first considering the model case $f=\chi_F$ and then elaborate on the
modifications needed to handle the general case.

\subsection{The case $f=\chi_F$}
$\newline$

Our intention here is to prove the following

\begin{t1}\label{lacchar} Let $F\subseteq\TT$ be a measurable set. Then
\beq\label{lacc}
\|T(\chi_F)\|_{L^{1,\infty}}\lesssim |F|\,\log\log(10+\frac{1}{|F|})\,\:.
\eeq
\end{t1}

With the notations from the previous section, statement \eqref{lacc} reduces to
\beq\label{lacc1}
\forall\:\:\:\bar{G}\subset\TT\:\:\:\:\:\:\:\:\:\:\:\:\|T^{*}(\chi_G)\|_{L^{1}(F)}\lesssim |F|\,\log\log(10+\frac{1}{|F|})\,\:.
\eeq

Thus, Theorem \ref{lacchar} will be a consequence of the following three propositions:

\begin{p1}\label{Cluster}
In the above settings we have
\beq\label{Clus}
\|\sum_{P\in\P_{cluster}} T_P^{*}(\chi_G)\|_{L^1(F)}\lesssim |F|\,\:.
\eeq
\end{p1}

\begin{p1}\label{Zygmund}
The following is true
\beq\label{Zyg}
\|\sum_{k\in\N}\sum_{I_{P_O}\in\I_k}\chi_{I_{P_O}}{T^{\p_k^{2}(P_O)}}^{*}(\chi_G)\|_{L^1(F)}\lesssim |F|\,\log\log(10+\frac{|G|}{|F|})\,\:.
\eeq
\end{p1}

\begin{p1}\label{p1fam}
The following relation holds
\beq \label{p1final}
\|\sum_{k\in\N}\sum_{I_{P_O}\in\I_k}\chi_{I_{P_O}}{T^{\p_{k}^{1}(P_O)}}^{*}(\chi_G)\|_{L^1(F)}\lesssim |F|\:.
\eeq
 \end{p1}

The proof of Proposition 1 is trivial if one uses the key observation that $\P_{cluster}$ is just an $\a-$dilation of a tree\footnote{See \cite{lv3} for definitions/notations.}. Indeed, heuristically, the proof reduces to the fact that the (maximal) Hilbert transform is bounded from $L^1$ to $L^{1,\infty}$. We leave the details for the reader.

For the second and third propositions, the strategy will be as follows: each $P_O$ generates a spacial band having the $x$ coordinate inside $I_{P_O}$. Thus given $P_O$, we will isolate the two corresponding families of tiles $\p_k^{2}(P_O)$ and $\p_k^{1}(P_O)$ respectively. For the family $\p_k^{2}(P_O)$, we restrict our analysis to the above mentioned band and apply orthogonality methods since we have enough separation - lacunary frequencies - among the tiles. For the tiles within $\p_k^{1}(P_O)$ we see no oscillation between the corresponding operators $\{\chi_{I_{P_O}}\,T_P^*\}_{P\in \p_k^{1}(P_O)}$ and thus we morally have ${T^{\p_k^{1}(P_O)}}^{*}(\chi_{G})(x)|_{x\in I_{P_O}}$ constant.

%Roughly, for $j\in\{1,2\}$, our aim is to get relations of the following type:
%\beq\label{ort}
% |\int_{F\cap I_{P_O}} {T^{\p_k^{j}(P_O)}}^{*}(\chi_{G})|\lesssim \l\,2^{-k}\,
% (k\,\log \frac{1}{\l})^{\frac{1}{2}}\,|I_{P_O}|^{\frac{1}{2}}\,\|\chi_{I_{P_O}}\,{T^{\p_k^{j}(P_O)}}^{*}(\chi_{G})\|_2^{\frac{1}{2}}\:.
%\eeq

 Before proceeding with the proof of our propositions we will need several notations.

Fix a collection $\p_k^{2}(P_O)$ and decompose it as a union of maximal trees $\bigcup_{l}\p_{k,l}(P_O)$ with each $\p_{k,l}(P_O)=\p_{k,l}$ a tree at the (dyadic) frequency $c_l$. Thus we have that
$${T^{\p_k^{2}(P_O)}}^{*}(\cdot)=\sum_{l} e^{-i\,c_l\cdot}\: {T^{\p_{k,l}^0}}^{*}(\cdot)\:,$$
where here $\p_{k,l}^0$ stands for the shift of $\p_{k,l}$ to the real axis.

For $I\subseteq\TT$ dyadic interval set
\beq\label{proj}
\L_{I}(f):=\frac{\int_{I}f(s)\,ds}{|I|}\chi_I\;.
\eeq

Define now the operator
\beq \label{almostconst}
{T_{c}^{\p_k^{2}(P_O)}}^{*}g(\cdot):=\sum_{l} e^{-i\,c_l\cdot}\: \L_{I_{P_O}}({T^{\p_{k,l}^0}}^{*}g)(\cdot)\:.
\eeq

\begin{l1}\label{treecutp}
The following holds
\beq \label{forestexpl2}
\|{T^{\p_k^{2}(P_O)}}^{*}(\chi_G)-{T_{c}^{\p_k^{2}(P_O)}}^{*}(\chi_G)\|_{L^{\infty} (I_{P_O})}\lesssim
\sum_{P\in\p_k^{2}(P_O)} \frac{|I_{P_O}|}{|I_P|}\frac{|E(P)\cap G|}{|I_P|}\;.
\eeq
 \end{l1}

\begin{proof}
%It will be enough to prove that there exists an absolute constant $C>0$ such that for any $p\in\N^{*}$ we have
%\beq \label{forestexpl2}
%\|\sum_{k} e^{-i\,c_k\cdot}\left({T^{\p_k^0}}^{*}\chi_{G}-\L_{\J}^{0}({T^{\p_k^0}}^{*}\chi_{G})\right)\|_{2p}\leq C\,\sqrt{p}\,2^{-\frac{n}{2}}\,\;.
%\eeq
Set for notational simplicity $T_l={T^{\p_{k,l}^0}}^{*}\chi_{G}-\L_{I_{P_O}}({T^{\p_{k,l}^0}}^{*}\chi_{G})$.

 For fixed $l$ and $x\in I_{P_O} $ we have
$$|T_l(x)|=\left|{T^{\p_{k,l}^0}}^{*}\chi_{G}(x)-\frac{1}{|I_{P_O}|}\int_{I_{P_O}}{T^{\p_{k,l}^0}}^{*}\chi_{G}(s)ds\right|=$$
$$\left|\frac{1}{|I_{P_O}|}\int_{I_{P_O}}\left\{\sum_{{P\in\p_{k,l}^0}\atop{2^{-j}=|I_P|\geq
|I_{P_O}|}}\int_{\TT}\left[\v_j(x-y)-\v_j(s-y)\right]\chi_{G}(y)\chi_{E(P)}(y)dy\right\}ds\right|$$
$$\lesssim\sum_{{P\in\p_{k,l}}\atop{I_{P^*}\supseteq I_{P_O}}} \frac{|I_{P_O}|}{|I_P|}\frac{|E(P)\cap G|}{|I_P|}\:.$$

Summing now in $l$ we deduce that \eqref{forestexpl2} holds.
\end{proof}

\begin{l1}\label{redconnst}
With the previous notations we have that
\beq \label{2reduction}
\|\sum_{k}\sum_{I_{P_O}\in\I_k}\chi_{I_{P_O}}{T^{\p_k^{2}(P_O)}}^{*}(\chi_G)\|_{L^1(F)}\lesssim \|\sum_{k}\sum_{I_{P_O}\in\I_k}{T_c^{\p_k^{2}(P_O)}}^{*}(\chi_G)\|_{L^1(F)}\,+\,|F|\:.
\eeq
 \end{l1}
\begin{proof}
Based on the previous lemma we deduce that
$$\|\sum_{k}\sum_{I_{P_O}\in\I_k}\chi_{I_{P_O}}{T^{\p_k^{2}(P_O)}}^{*}(\chi_G)\|_{L^1(F)}\lesssim $$ $$\|\sum_{k}\sum_{I_{P_O}\in\I_k}{T_c^{\p_k^{2}(P_O)}}^{*}(\chi_G)\|_{L^1(F)}\,+\,
\:\sum_{k}2^{-k}\,\l\,\sum_{I_{P_O}\in\I_k}\sum_{{P\in\bigcup_{l}\p_{k,l}(P_O)}} \frac{|I_{P_O}|^2}{|I_P|^2}|E(P)\cap G|.$$
Thus, for proving \eqref{2reduction} it will be enough to prove that
\beq \label{kdiscard}
2^{-k}\,\l\,S_{k}\lesssim 2^{-k/2}\,|F|\:,
\eeq
where
$$S_k:=\sum_{I_{P_O}\in\I_k}\sum_{{P\in\bigcup_{l}\p_{k,l}(P_O)}} \frac{|I_{P_O}|^2}{|I_P|^2}|E(P)\cap G|=\sum_{J\in\I_{k+1}}S_{k,J}\:,$$
and for a fixed $J\in\I_{k+1}$ we set
$$S_{k,J}:=\sum_{{I_{P_O}\in\I_k}\atop{I_{P_O}\subset J}}\sum_{{P\in\bigcup_{l}\p_{k,l}(P_O)}\atop{J\supseteq I_{P^*}^r\supseteq I_{P_O}}}
\frac{|I_{P_O}|^2}{|I_P|^2}|E(P)\cap G|\:.$$
Further, defining
$$\A_{m, J}:=\{I\:\textrm{dyadic}\,|\,100J\supseteq 20I\:\:\textrm{and}\:\:\frac{|I\cap G|}{|I|}\approx 2^{-m}\}\:\:\:\textrm{with}\:\:\:m\in \N\:,$$
set
$$S_{k,J}^{m}:=\sum_{{I_{P_O}\in\I_k}\atop{I_{P_O}\subset J}}\sum_{{P\in\bigcup_{l}\p_{k,l}(P_O)}\atop{J\supseteq I_{P^*}^r\supseteq I_{P_O},
\:I_{P}\in \A_{m, J}}}\frac{|I_{P_O}|^2}{|I_P|^2}|E(P)\cap G|\:,$$
and respectively
$$L_{k,J}^{m}:=\sum_{{I_{P_O}\in\I_k}\atop{I_{P_O}\subset J}}\sum_{{20I\supseteq I_{P_O}}\atop{I\in \A_{m, J}}}
\frac{|I_{P_O}|^2}{|I|^2}|I\cap G|\:.$$
Remark that
$$S_{k,J}^{m}\leq L_{k,J}^{m}\,$$
and thus if
$$L_{k,J}:=\sum_{{I_{P_O}\in\I_k}\atop{I_{P_O}\subset J}}\sum_{{100J\supseteq 20I\supseteq I_{P_O}}\atop{I\:\textrm{dyadic}}}
\frac{|I_{P_O}|^2}{|I|^2}|I\cap G|\:,$$
then we have that
$$S_{k,J}=\sum_{m\in\N}S_{k,J}^{m}\lesssim L_{k,J}=\sum_{m\in \N}L_{k,J}^{m}\;.$$
Now set $\A_{m, J}^{max}:=\{I\in \A_{m, J}\,|\,I \textrm{maximal}\}$ and deduce that
$$L_{k,J}^{m}\lesssim 2^{-m}\,\sum_{{I_{P_O}\in\I_k}\atop{I_{P_O}\subset J}}\sum_{{20I\supseteq I_{P_O}}\atop{I\in \A_{m, J}}}
\frac{|I_{P_O}|}{|I|}|I_{P_0}|\lesssim  2^{-m}\,\sum_{{I_{P_O}\in\I_k}\atop{I_{P_O}\subset \bigcup_{I\in \A_{m, J}^{max}}20 I}}\,|I_{P_0}|$$
$$\lesssim 2^{-m}\,\sum_{I\in \A_{m, J}^{max}} |I|\lesssim 2^{-m/2}\,\sum_{I\in \A_{m, J}^{max}} |I|^{\frac{1}{2}}\,|I\cap G|^{\frac{1}{2}}\:,$$
from which we conclude
\beq \label{am}
L_{k,J}^{m}\lesssim 2^{-m/2}\,|J|^{\frac{1}{2}}\, |100 J\cap G|^{\frac{1}{2}}\:.
\eeq
It will be  thus enough to show that\footnote{Here $100\,\bar{\I}_{k+1}:=\bigcup_{J\in\I_{k+1}} 100\,J$.}
\beq \label{msum}
\sum_{J\in\I_{k+1}}|J|^{\frac{1}{2}}\,|100\,J\cap G|^{\frac{1}{2}}\lesssim
|\bar{\I}_{k+1}|^{\frac{1}{2}}\,|100\,\bar{\I}_{k+1}\cap G|^{\frac{1}{2}}
\eeq

 For this we will apply a greedy algorithm that may be regarded as an iteration of a Vitali covering type argument:

 Select the largest\footnote{If there are two (or more) intervals with maximal length just choose one of them.} interval $J_1$ inside the set $\I_{k+1}$, take its enlargement $100 J_1$, and let $\B_{1,1}$ be the set of all intervals $J\in\I_{k+1}$ such that $100 J\cap 100 J_1\not=\emptyset$. Then, repeat this procedure for the set  $\I_{k+1}\setminus \B_{1,1}$ thus obtaining a (maximal) interval $J_2$ and a set of subordinate intervals $\B_{1,2}$. Take now the set $\I_{k+1}\setminus (\B_{1,1}\cup \B_{1,2})$ and repeat this procedure till exhaustion.\footnote{Wlog we may assume that the set $\I_{k+1}$ is finite and hence this selection algorithm will finish in a finite number of steps.} Denote with $\B_1$ the set of the maximal intervals $\{J_l\}_{l=1}^{N}$ obtained at the moment of our stopping time (exhaustion). With this done, define $\I_{k+1}^1:=\I_{k+1}\setminus \B_1$.

 Repeat the entire algorithm described above for this set of intervals and obtain a new defined set $\B_2$. Let $\I_{k+1}^2:=\I_{k+1}\setminus (\B_1\cup\B_2)$ and apply again this algorithm. This procedure will end up in $p\in\N$ steps.

 In this way we were able to obtain a partition of the set
 $$\I_{k+1}=\bigcup_{r=1}^{p} \B_r\:,$$
such that for each $r\in\{1,\ldots,p\}$ the set $\{100 J\}_{J\in\B_r}$ consists of disjoint intervals.

Moreover, denoting with $\bar{\B}_r:=\bigcup_{J\in\B_r}J$, we have the following key relation:
\beq \label{ineq}
|\bar{\B}_l|\geq\frac{1}{500} \sum_{r=l}^{p} |\bar{\B}_r|\:\:\:\:\:\:\textrm{for any}\:\:\:l\in\{1,\ldots, p\}\:.
\eeq
Indeed, we have that
$$|\bar{\B}_l|=\frac{1}{500}\sum_{J\in\B_l}5\,|100J|\geq\frac{1}{500}\sum_{J\in\B_l}\sum_{{100 I\cap 100 J\not=\emptyset}\atop{I\in\I_{k+1}\setminus \bigcup_{r=1}^{l-1}\B_r}}|I|=\frac{1}{500}\sum_{r=l}^{p} |\bar{\B}_r|\:.$$
Now, from \eqref{ineq}, for any $s\in\N$ such that $p\geq 1000 (s+1)$ we have
\beq \label{bound}
|\bar{\B}_{1000s+1}|+\cdot+|\bar{\B}_{1000s+1000}|\leq\frac{1}{2^s} (|\bar{\B}_{1}|+\cdot+|\bar{\B}_{1000}|)\:.
\eeq

Assume wlog that $p=1000 s_0$ for some $s_0\in\N$. From \eqref{bound} and Cauchy-Schwarz inequality we conclude
$$\sum_{J\in\I_{k+1}}|J|^{\frac{1}{2}}\,|100J\cap G|^{\frac{1}{2}}\leq\sum_{s=0}^{s_0-1} \sum_{r=1000s+1}^{1000(s+1)}
 |\bar{\B}_r|^{\frac{1}{2}}\,|100\I_{k+1}\cap G|^{\frac{1}{2}}$$
 $$\lesssim |100\I_{k+1}\cap G|^{\frac{1}{2}}\,\sum_{s=0}^{s_0-1}[\frac{1}{2^s} (|\bar{\B}_{1}|+\cdot+|\bar{\B}_{1000}|)]^{\frac{1}{2}}
 \lesssim |\I_{k+1}|^{\frac{1}{2}}\,|100\I_{k+1}\cap G|^{\frac{1}{2}}\;.$$
This ends the proof of \eqref{msum}.
\end{proof}

Now let $\p_k^{2}(P_O)=\bigcup_{n\in\N}\p_{k,n}^{2}(P_O)$ be the mass decomposition of $\p_k^{2}(P_O)$ referred to in Section 3.1..
As before, we decompose $\p_{k,n}^{2}(P_O)$ in a union of maximal trees $\bigcup_{l}\p_{k,n,l}$ with each $\p_{k,n,l}$ a tree at the (dyadic) frequency $c_l$. Further, set the square function
\beq \label{sqf}
S_{\p_{k,n}^{2}(P_O)}:=\{\sum_{l}|{T^{\p_{k,n,l}}}^{*}|^2\}^{\frac{1}{2}}\:.
\eeq
Then, we have the following:

\begin{l1}\label{kl2zyg}
For each $n\in\N$ we have that
\beq \label{k2orthog}
\|{T_c^{\p_{k,n}^{2}(P_O)}}^{*}(\chi_G)\|_{L^1(F)}\lesssim 2^{-k}\,\l\,(k\log\frac{1}{\l})^{\frac{1}{2}}\,
|I_{P_O}|^{\frac{1}{2}}\,\|S_{\p_{k,n}^{2}(P_O)}(\chi_G)\|_{L^2(I_{P_O})}\:.
\eeq
 \end{l1}
\begin{proof}
The proof of this lemma will be based on the following facts:
\begin{itemize}
\item the \textit{good separation} among the frequencies of the trees in $\p_{k,n}^{2}(P_O)$ - this is a consequence of the $(f,\l)$-lacunary decomposition described in the previous section;
\item \textit{Zygmund's inequality} regarding the behavior of the lacunary Fourier series with $l^2-$coefficients:
\beq \label{zyg}
\|\sum_{j} a_j\,e^{i\,n_j\,x}\|_{\exp(L^2(\TT))}\lesssim \{\sum_{j} |a_j|^2\}^{\frac{1}{2}}\,,
\eeq
where here $\{n_j\}_j\subset\N$ is a lacunary sequence.
\end{itemize}

We mention here that the idea of using Zygmund's inequality in the context of Konyagin's question (for the Walsh model) was appearing in \cite{LaDo}. As it turns out, while beautiful in nature, this ingredient is not actually required\footnote{For more details see Remarks section.} for obtaining our theorem, though it offers the best bound in \eqref{k2orthog}.

Let us pass now to the actual proof of our lemma. With the previous notations we set
$${T_{c}^{\p_{k,n}^{2}(P_O)}}^{*}(\chi_G)(\cdot):=\chi_{I_{P_O}}(\cdot)\,\sum_{l} e^{-i\,c_l\cdot}\, \L_{I_{P_O}}({T^{\p_{k,n,l}^0}}^{*}(\chi_G))(\cdot)\:.$$
Thus applying twice Cauchy-Schwartz inequality followed by (the dual form of) Zygmund inequality we conclude
$$\left|\int_{F\cap I_{P_O}}\sum_{l} e^{-i\,c_l\cdot}\, \L_{I_{P_O}}({T^{\p_{k,n,l}^0}}^{*}(\chi_G))\right|$$
$$\lesssim\|\{\sum_{l}|{T^{\p_{k,n,l}^0}}^{*}(\chi_G)|^2\}^{\frac{1}{2}}
\{\sum_l|\L_{I_{P_O}}(e^{-i\,c_l\cdot}\,\chi_F)|^2\}^{\frac{1}{2}}\|_{L^1(I_{P_O})}$$
$$\lesssim \|S_{\p_{k,n}^{2}(P_O)}(\chi_G)\|_{L^2(I_{P_O})}\,\left\{\sum_l\frac{|<\chi_{F\cap I_{P_O}},\,e^{-i\,c_l\cdot} >|^2}{|I_{{P_O}}|}\right\}^{\frac{1}{2}}$$
$$\lesssim\frac{|I_{P_O}\cap F|}{|I_{P_O}|}\,\left(\log \frac{|I_{P_O}\cap F|}{|I_{P_O}|}\right)^{\frac{1}{2}}\,
|I_{P_O}|^{\frac{1}{2}}\,\|S_{\p_{k,n}^{2}(P_O)}(\chi_G)\|_{L^2(I_{P_O})}\;,$$
which proves \eqref{2orthog}.
\end{proof}

\begin{l1}\label{allkl2zyg}
With the previous notations we have
\beq \label{2orthog}
\|\sum_{k\in\N}\sum_{I_{P_O}\in\I_k}{T_c^{\p_{k,n}^{2}(P_O)}}^{*}(\chi_G)\|_{L^1(F)}\lesssim 2^{-n/2}\,|F|\,(\log\frac{|G|}{|F|})^{\frac{1}{2}}\,\:.
\eeq
 \end{l1}
\begin{proof}
This proof is a consequence of relation \eqref{k2orthog} and of the $L^2-$forest estimate (see e.g. main Theorem, c), in \cite{lv7}). Indeed, based on these, we have that
 $$\sum_{I_{P_O}\in\I_k}\|{T_c^{\p_{k,n}^{2}(P_O)}}^{*}(\chi_G)\|_{L^1(F)}$$
 $$\lesssim \sum_{I_{P_O}\in\I_k} 2^{-k}\,\l\,(k\log\frac{1}{\l})^{\frac{1}{2}}\,|I_{P_O}|^{\frac{1}{2}}\,\|S_{\p_{k,n}^{2}(P_O)}(\chi_G)\|_{L^2(I_{P_O})}$$
 $$\leq2^{-k}\,\l\,(k\log\frac{1}{\l})^{\frac{1}{2}}\,|\bar{\I}_k|^{\frac{1}{2}}\,
 \|\sum_{I_{P_O}\in\I_k}\chi_{I_{P_O}}\,S_{\p_{k,n}^{2}(P_O)}(\chi_G)\|_{L^2}$$
 $$\lesssim 2^{-n/2}\,(2^{-k}\,\l)^{\frac{1}{2}}\,(k\log\frac{1}{\l})^{\frac{1}{2}}\,|2^{-k}\,\l\,\bar{\I}_k|^{\frac{1}{2}}\,|G|^{\frac{1}{2}}
 \lesssim 2^{-n/2}\,(2^{-k}\,k)^{\frac{1}{2}}\,(\log\frac{1}{\l})^{\frac{1}{2}}\,|F|\:.$$
 Summing now in $k\in\N$ and using triangle inequality we deduce that \eqref{2orthog} holds.
\end{proof}

The next lemma follows from inspecting the proof of part b) of the main theorem in \cite{lv7}:

\begin{l1}\label{l1weak} The following holds

\beq \label{l1w}
\|\sum_{k\in\N}\sum_{I_{P_O}\in\I_k}{T_c^{\p_{k,n}^{2}(P_O)}}^{*}(\chi_G)\|_{L^1(F)}\lesssim |F|\:.
\eeq
 \end{l1}
$\newline$

\noindent\textbf{Proof of Proposition 2.}
$\newline$

In the view of Lemmas \ref{redconnst}, \ref{allkl2zyg} and \ref{l1weak} we deduce that

$$\|\sum_{k\in\N}\sum_{I_{P_O}\in\I_k}\chi_{I_{P_O}}{T^{\p_k^{2}(P_O)}}^{*}(\chi_G)\|_{L^1(F)}$$
$$\lesssim |F|\,+\,|F|\,\sum_{n\in\N}\min\{2^{-n/2}\,(\log\frac{|G|}{|F|})^{\frac{1}{2}},\,1\}\lesssim |F|\,\log\log(10+\frac{|G|}{|F|})\:.$$
Thus \eqref{Zygmund} holds.
$\newline$
$\newline$

We pass now to the proof of Proposition 3.

Here we are using the simple observation that on the interval $I_{P_O}$ our operator ${T^{\p_k^{1}(P_O)}}^{*}(\chi_{G})$ is morally constant.

For showing this we first proceed as in the case of $\p_k^{2}(P_O)$, and decompose $\p_k^{1}(P_O)$ as a union of maximal trees $\bigcup_{l}\p_{k,l}(P_O)$ with each $\p_{k,l}(P_O)=\p_{k,l}$ a tree at the (dyadic) frequency $c_l$. Here one should notice the key fact that from the definition of $\p_k^{1}(P_O)$ we may wlog suppose that for all $l$ in the above decomposition we have
\beq \label{smallosc}
c_l\,|I_{P_O}|\leq\frac{1}{2}\;.
\eeq
Now as before, set
$${T^{\p_k^{1}(P_O)}}^{*}=\sum_{l} e^{-i\,c_l\cdot} {T^{\p_{k,l}^0}}^{*}\:.$$
and define
\beq \label{almostconst1}
{T_{c}^{\p_k^{1}(P_O)}}^{*}g(\cdot):=\chi_{I_{P_O}}(\cdot)\,\sum_{l} e^{-i\,c_l\cdot}\, \L_{I_{P_O}}({T^{\p_{k,l}^0}}^{*}g)(\cdot)\:.
\eeq
Now following the same steps from the treatment of $\p_k^{2}(P_O)$ we deduce that
\begin{itemize}
\item  as in Lemma \ref{treecutp} we have

\beq \label{forestexpl1}
\|{T^{\p_k^{1}(P_O)}}^{*}(\chi_G)-{T_{c}^{\p_k^{1}(P_O)}}^{*}(\chi_G)\|_{L^{\infty} (I_{P_O})}\lesssim
\sum_{{P\in\p_k^{2}(P_O)}} \frac{|I_{P_O}|}{|I_P|}\frac{|E(P)\cap G|}{|I_P|}\;.
\eeq

\item as in Lemma \ref{redconnst} we have

\beq \label{1reduction}
\sum_{k}\|\sum_{I_{P_O}\in\I_k}\chi_{I_{P_O}}{T^{\p_k^{1}(P_O)}}^{*}(\chi_G)\|_{L^1(F)}\lesssim \sum_{k}\|\sum_{I_{P_O}\in\I_k}{T_c^{\p_k^{1}(P_O)}}^{*}(\chi_G)\|_{L^1(F)}\,+\,|F|\:.
\eeq

\end{itemize}

Using now \eqref{smallosc} together with the Taylor expansion
$$e^{i\,c_l\,(x-y_{P_O})}:=\sum_{k\geq 0} \frac{[c_l\,(x-y_{P_O})]^k}{k!}$$ which
is absolutely and uniformly convergent for any $x,\,y_{P_O}\in I_{P_O}$ we deduce that
\beq \label{const1}
\begin{array}{rl}
\|{T_c^{\p_k^{1}(P_O)}}^{*}(\chi_G)\|_{L^1(F)}&\lesssim |I_{P_O}\cap F|\,|\sum_{l} e^{-i\,c_l\,y_{P_O}}\, \L_{I_{P_O}}({T^{\p_{k,l}^0}}^{*}(\chi_G))(y_{P_0})|\\
&+\,|I_{P_O}\cap F|\,\sup_{l}\frac{\int_{I_{P_O}}|{{T^{\p_{k,l}}}^{*}(\chi_G)}|}{|I_{P_O}|}\;.
\end{array}
\eeq
Thus, from \eqref{1reduction} and \eqref{const1}, for an appropriately chosen $y_{P_O}\in I_{P_O}$ we deduce
\beq \label{const12}
\begin{array}{cl}
\sum_{k}\|\sum_{I_{P_O}\in\I_k}{T_c^{\p_k^{1}(P_O)}}^{*}(\chi_G)\|_{L^1(F)}\lesssim \sum_{k} 2^{-k}\,\l\,\|\sum_{I_{P_O}\in\I_k}{T_c^{\p_k^{1}(P_O)}}^{*}(\chi_G)\|_{L^1}\\
+\,\sum_{k} 2^{-k}\,\l\,\sum_{I_{P_O}\in\I_k}\sup_{l}\int_{I_{P_O}}|{{T^{\p_{k,l}(P_O)}}^{*}(\chi_G)}|\:.
\end{array}
\eeq
Now using \eqref{const12} and applying \eqref{forestexpl1}, \eqref{1reduction} backwards we deduce
$\newline$

\begin{l1}\label{firstaproxp1fam}
The following holds
\beq \label{1noosc}
\begin{array}{rl}
\|\sum_{k\in\N}\sum_{I_{P_O}\in\I_k}&{T^{\p_{k}^{1}(P_O)}}^{*}(\chi_G)\|_{L^1(F\cap I_{P_O})}\\
&\lesssim\sum_{k} 2^{-k}\,\l\,\|\sum_{I_{P_O}\in\I_k}{\chi_{I_{P_O}}\,T^{\p_k^{1}(P_O)}}^{*}(\chi_G)\|_{L^1}\\
&+\,\sum_{k} 2^{-k}\,\l\,\sum_{I_{P_O}\in\I_k}\sup_{l}\int_{I_{P_O}}|{{T^{\p_{k,l}(P_O)}}^{*}(\chi_G)}|\,+\,|F|\:.

\end{array}
\eeq
 \end{l1}
$\newline$
$\newline$

\noindent\textbf{Proof of Proposition 3.}
$\newline$

Applying now the mass decomposition (see Section 3.1.) $\p_k^{1}(P_O)=\bigcup_{n\in\N}\p_{k,n}^{1}(P_O)$, \textit{i.e.} the decomposition of $\p_k^{1}(P_O)$ into families with uniform mass parameter, we observe that a similar relation with \eqref{1noosc} will hold for the corresponding $n-$level.

Using the $L^2$-bound on each $n-$th generation forest (see e.g. \cite{lv7}) we have
$$ 2^{-k}\,\l\,\left(\|\sum_{I_{P_O}\in\I_k}{\chi_{I_{P_O}}\,T^{\p_{k,n}^{1}(P_O)}}^{*}(\chi_G)\|_{L^1}
+\sum_{I_{P_O}\in\I_k}\sup_{l}\int_{I_{P_O}}|{{T^{\p_{k,n,l}(P_O)}}^{*}(\chi_G)}|\right) $$
$$\lesssim 2^{-k/2}\,2^{-n/2}\,(2^{-k}\,\l\,|\I_k|)^{\frac{1}{2}}\,(\l\,|G|)^{\frac{1}{2}}\lesssim 2^{-k/2}\,2^{-n/2}\,|F|\,,$$
and thus we conclude that \eqref{p1final} holds.
$\newline$

\subsection{The general case}
$\newline$

Let us now suppose that $f\in L^1(\TT)$. We will then modify accordingly the $(f,\l)-$decomposition presented in Section 3.2.. Beyond the standard adaptations we only need to notice that while formulation \eqref{tot} is preserved\footnote{Again the reader should notice that this decomposition
of the family $\P$ does \textit{not} form a partition.}
$$\P=\P_0\cup\P_{cluster}\cup\bigcup_{k\in\N}\bigcup_{P_O\in\p_{k,O}}\p_k^{2}(P_O)\cup\p_k^{1}(P_O)\:,$$
the definition of $\P_0$ in \eqref{rest} should be now rephrased as follows:
\beq\label{restref}
\P_0=\bigcup_{l\in\N}\P_0(l)\;,
\eeq
where
$$\P_0(l):=\{P\in\P\,|\,\exists\:I_{P^*}^r\subseteq I\in\I_l\:\:\&\:\:I_{P^*}^r\cap\bar{\I}_{l-1}=\emptyset\}\:\:\:\textrm{for}\:l\geq 1\:,$$
and
$$\P_0(0):=\left\{P\in\P\,|\,I_{P^*}\subset 1000\,\{x\,|\,M(f)(x)> \frac{1}{2}\,\l\}\right\}\:.$$

 Next, let us notice the simple but key observation that for any $I$ (dyadic) interval with $I\in\bar{\I}_l$ and $I\cap\bar{\I}_{l-1}=\emptyset$ (here $l\geq 1$) one has the pointwise estimate
$$|f(x)|\leq 2^{-l+1}\,\l\:\:\:\textrm{a.e.}\:\:x\in I\;.$$
This way we will be able for tiles in $\P_0(l)$ to obtain the desired estimates as a consequence of the $L^{\infty}-$bounds of $f$.

More precisely, the analogue of \eqref{partition} becomes now

\beq\label{partitiongen}
\begin{array}{rl}
|\int f\,\sum_{P\in\P} T_P^{*}(\chi_G)|&\lesssim |\int f\,\sum_{P\in\P_{cluster}} T_P^{*}(\chi_G)|\\
&+\,|\int f\,\sum_{k}\sum_{I_{P_O}\in\I_k}\chi_{I_{P_O}}{T^{\p_k^{2}(P_O)}}^{*}(\chi_G)|\\
&+\,|\int f\,\sum_{k}\sum_{I_{P_O}\in\I_k}\chi_{I_{P_O}}{T^{\p_k^{1}(P_O)}}^{*}(\chi_G)|\\
&+\,\sum_{k\in\N}\,2^{-k}\,\l\,\int_{\bar{\I}_k}|T_k^{*}(\chi_G)|\;,
\end{array}
\eeq
where here we set $T_k^{*}:=\sum_{P\in\P} \chi_{P}^k\,T_P^{*}$.

In the above definition we impose $\chi_{P}^k(x)=1$ if $\exists\,r$ such that $x\in I_{P^*}^r\subseteq\bar{\I}_k$ and $I_{P^*}^r\cap\bar{\I}_{k-1}=\emptyset$; otherwise we set $\chi_{P}^k(x)=0$.

Using the $L^2$ boundedness of the Carleson operator for the fourth term in the right hand term summation we obtain
the upper bound estimate
$$\sum_{k\in\N} 2^{-k}\,\l\,|\bar{\I}_k|^{\frac{1}{2}}\,|100\bar{\I}_k\cap G|^{\frac{1}{2}}\lesssim \sum_{k}2^{-\frac{k}{2}}\,\|f\|_1\lesssim \|f\|_1\:.$$

With these being said, one can modify the above proofs of Propositions 1, 2 and 3 accordingly:
$\newline$

\noindent\textbf{Proposition 1$^{\prime}$.}$\:$\textit{With the above notations we have}
\beq\label{Clusx}
\left|\int f\,\sum_{P\in\P_{cluster}} T_P^{*}(\chi_G)\right|\lesssim \|f\|_1\,\:.
\eeq

\noindent\textbf{Proposition 2$^{\prime}$.}$\:$\textit{The following is true}
\beq\label{Zygx}
\left|\int f\,\sum_{k\in\N}\sum_{I_{P_O}\in\I_k}\chi_{I_{P_O}}{T^{\p_k^{2}(P_O)}}^{*}(\chi_G)\right|\lesssim \|f\|_1\,\log\log(10+\frac{|G|\|f\|_{\infty}}{\|f\|_1})\,\:.
\eeq

\noindent\textbf{Proposition 3$^{\prime}$.}$\:$\textit{The following relation holds}
\beq \label{p1ffinal}
\left|\int f\,\sum_{k\in\N}\sum_{I_{P_O}\in\I_k}\chi_{I_{P_O}}{T^{\p_{k}^{1}(P_O)}}^{*}(\chi_G)\right|\lesssim \|f\|_1\:.
\eeq

While Propositions $1^{\prime}$ and $3^{\prime}$ only require trivial modifications for Proposition $2^{\prime}$ one needs to modify the proof of Lemma \ref{allkl2zyg} and show that
 $$\sum_{I_{P_O}\in\I_k}\left|\int f\,{T_c^{\p_{k,n}^{2}(P_O)}}^{*}(\chi_G)\right| \lesssim 2^{-n/2}\,2^{-k/2}\,(\log\frac{\|f\|_{\infty}}{2^{-k}\,\l})^{\frac{1}{2}}\,\|f\|_1\:.$$
This is a direct consequence of the fact that for $I_{P_O}\in\I_k$ and $\textbf{c}=\{c_l\}_l\subset \N$ lacunary one has
$$\left\{\sum_l\frac{|\int_{I_{P_O}}f\,e^{-i\,c_l\cdot} |^2}{|I_{{P_O}}|}\right\}^{\frac{1}{2}}\lesssim_{\textbf{c}} 2^{-k}\,\l\,\left(\log\frac{\|f\|_{\infty}}{2^{-k}\,\l}\right)^{\frac{1}{2}}\,|I_{P_O}|^{\frac{1}{2}}\:.$$
We leave further details for the interested reader.

With these done we deduce that the following holds:

\begin{t1}\label{lacatom}\textit{ Let $f\in L^{1}(\TT)\cap L^{\infty}(\TT)$. Then we have}
\beq\label{lacat}
\|T(f)\|_{L^{1,\infty}}\lesssim \|f\|_1\,\log\log(10+\frac{\|f\|_{\infty}}{\|f\|_{1}})\,\:.
\eeq
\end{t1}

This further implies our main result:

\begin{t1}\label{lacconv} Let $f\in L\log\log L\log\log\log L (\TT)$. Then we have
\beq\label{lac}
\|T(f)\|_{L^{1,\infty}}\lesssim \|f\|_{L\log\log L\log\log\log L }\,\:.
\eeq
\end{t1}

\begin{proof}

Let $\W$ be the quasi-Banach space defined as follows:
$$\W:=\{f:\:\TT\mapsto C\,|\,f\:\textrm{measurable},\:\|f\|_{\W}<\infty \}\:\:\:\textrm{where}$$
$$\|f\|_{\W}:=\inf\left\{\sum_{j=1}^{\infty}(1+\log j)\|f_j\|_1\,\log\log\frac{e\,\|f_j\|_{\infty}}{\|f_j\|_1}\:\:\left|\right.\:\:
\begin{array}{cl}
f=\sum_{j=1}^{\infty}f_j,\:\\
\sum_{j=1}^{\infty}|f_j|<\infty\:\textrm{a.e.}\\
f_j\in L^{\infty}(\TT)
\end{array}
\right\}\;.$$

Next, using a similar reasoning as in \cite{Ar}, one can show that the following holds:
\beq\label{w}
\|\cdot\|_{\W}\leq\|\cdot\|_{L\log\log L\log\log\log L }\:.
\eeq

Let now $f=\sum_{j=1}^{\infty}f_j$ a decomposition of $f$ as described in the $\W$-norm definition. Then applying Kalton's inequality (\cite{Ka})
and Theorem \ref{lacatom} we have
$$\|T(f)\|_{L^{1,\infty}}\lesssim\sum_{j\geq 1}(1+\log j)\,\|T(f_j)\|_{L^{1,\infty}}$$
$$\lesssim \sum_{j\geq 1}(1+\log j)\,\|f_j\|_1\,\log\log(10+\frac{\|f_j\|_{\infty}}{\|f_j\|_{1}})\,,$$
which based on \eqref{w} proves \eqref{lac}.
\end{proof}

\section{Remarks}

1) For approaching Konyagin's conjecture, in both \cite{LaDo} (for the Walsh case) and our paper (for the Fourier case), one follows several natural steps given the formulation of the problem: the space $L\log \log L$ is regarded as an intermediate space between the space $L^{1,\infty}$ and some space of the form $L(\log L)^{\a}$ with $\a>0$.

The $L^{1,\infty}$ space appears when bounding the tile-families with uniform mass (and size) and its usage is suggested by the $L^1$-weak boundedness of the Hilbert transform encoded in $H: \:L^1\rightarrow L^{1,\infty}$.

The other ``end-point" of the spectrum - the space $L(\log L)^{\a}$ - appears as a manifestation of the lacunary structure of the maximal operator $S_{lac}$ and is the key place where one makes use of the specific nature of our problem. The heart of the matter here, is represented by the fact
that the sequence of torus characters $\{e^{i\,2^j\,x}\}_{j\in\N}$ behaves as good as a sequence of i.i.d. random variables. Indeed, if $(\Omega,\f, \mu )$ is a probability space and $\{r_j\}_j$ is a sequence of independent random variables taking values $\pm 1$ with equal probability then
Khinchin's inequality  asserts that for any $0<p<\infty$ and $\{a_j\}_j\in l^2(\N)$
\beq\label{Kh1}
\|\sum_{j} a_j\,r_j\|_{L^p(\Omega)}\lesssim_{p}\|\{a_j\}_j\|_{l^2(\N)}\:.
\eeq
 Moreover, for some $c>0$, one has\footnote{See \textit{e.g.} \cite{W}.}
\beq\label{Kh2}
\mu(\{t\,|\,|\sum_{j} a_j\,r_j(t)|>\l\})\lesssim e^{-\frac{c\,\l^2}{\|\{a_j\}_j\|^2_{l^2(\N)}}} \:.
\eeq
Both \eqref{Kh1} and \eqref{Kh2} remain valid when replacing $\{r_j\}_j,\,\Omega,\mu$ with the counterparts $\{e^{i\,2^j\,x}\}_j,\,[0,1],\,dx$!\footnote{The proof of
this fact is a real analysis exercise and reduces to showing that $\|\sum_{j} a_j\,e^{i\,2^j\,\cdot}\|_{L^p(\TT)}\lesssim p^{\frac{1}{2}}\,\|\{a_j\}_j\|_{l^2(\N)}$ for $p\in2\N$.} Now the Fourier
version of \eqref{Kh2} is precisely the Zygmund inequality \eqref{zyg}, which viewed dually may be written as\footnote{Here $\hat{f}(n)$ stands for the $n^{th}$ Fourier coefficient of $f$.}
\beq\label{Kh2dual}
\|\{\hat{f}(2^j)\}_j\|_{l^2(\N)}\lesssim \|f\|_{L\,(\log L)^{\frac{1}{2}}(\TT)}\:,
\eeq
thus reaching the announced target space.

This upgrade of the Hausdorff-Young inequality is the point that makes possible an estimate like \eqref{K1}, far below the critical $L\log L$ space as appearing in Conjecture 1.

As a last remark on this theme, one should notice that the statement
\beq\label{Kh3dual}
\|\{\hat{f}(2^j)\}_j\|_{l^2(\N)}\lesssim \|f\|_{L\,(\log L)^{\a}(\TT)}\:\:\:\textrm{for some}\:\a>0,
\eeq
is still enough for proving our theorem. The fact that \eqref{Kh3dual} holds for $\a\geq1$ is a simple consequence of the
relation\footnote{See \textit{e.g.} \cite{S}, p. 178.}
\beq\label{BMO}
\{a_j\}_j\in l^2(\N)\:\:\:\Rightarrow\:\:\:f(x):=\sum_{j\in\N} a_j\,e^{i\,2^j\,x}\in BMO(\TT)\:.
\eeq

2) Besides Zygmund's inequality, the other ingredients used in \cite{LaDo} involve the time-frequency
approach developed in \cite{lt3}, a multi-frequency projection argument in the spirit of \cite{NOT}, rearrangement invariant spaces techniques and extrapolation theory gradually built on the results in \cite{sj2}, \cite{So2}, \cite{An}, \cite{SS}, and \cite{CM}.

In our paper however, we take a different path with the tile decomposition in Section 3 playing the central role in the mechanism of our proof. The methods developed in \cite{lv7} are also relevant. Below, we present part of an antithesis between \cite{LaDo} and the present paper:
\begin{itemize}
\item \textit{in \cite{LaDo}, the authors are making an essential use of the following property: if $w_P$ and $w_{P_1}$ are Walsh wave packets adapted
to the tiles $P$, $P_1$ then the condition  $P<P_1$ implies $w_{P_1}=C\,w_{P}$ on $I_{P_1}$ where $C$ is a constant. Indeed, this property facilitates an application of a projection argument which does not seem to have a direct correspondent in the continuous case.}

This property greatly simplifies the picture in the Walsh setting and, in an artificial manner from the Fourier setting point of view, gets rid of the family of tiles which are not well separated within the scale imposed by the exceptional set. This is precisely why one needs to introduce the $(f,\l)$ decomposition presented in Section 3 and further why one needs to treat separately the cases described by Propositions 1, 2 and 3.
\item \textit{in our paper we embrace the approach introduced by Fefferman in \cite{f}, and further developed in \cite{lv3}, \cite{lv7}}.

Among others, this offers us more flexibility in treating the forest estimates and the advantage of having a tile decomposition which is directly adapted to the exceptional sets of the Carleson operator and to the size of $f$.
\item \textit{based on the techniques developed in \cite{lv7}, we are able to eliminate the use of extrapolation theory in passing from restricted weak type to just weak type estimates}.
\end{itemize}

3) The remaining gap between our result - the pointwise convergence of the lacunary Fourier series in $L\log\log L\log\log\log L$ - and the conjectured (best possible) space $L\log\log L$ is generated by the same\footnote{One may want to compare the parallelism between the approach presented here and the one in paper \cite{lv7}.} (lack of) technology as in the case of the pointwise convergence of the full sequence of the partial Fourier series - \textit{e.g.} Antonov's result on the convergence in $L\log L\log\log\log L$ versus the conjectural space $L\log L$. Thus any progress on one of the problems will very likely imply a similar progress on the other.
$\newline$

\noindent\textbf{Acknowledgements}: I would like to thank Yen Do for a useful conversation on his joint work with Michael Lacey on the convergence of the lacunary Walsh-Fourier Series that motivated me further in finding a solution to the problem exposed in the present paper. Also, we thank Michael Lacey for some useful comments on an earlier draft of this paper.

\end{document}